\documentclass[twoside,12pt,leqno]{article}
\advance\topmargin by -\headheight\advance\topmargin by -\headsep

\usepackage{amsmath}
\usepackage{amsthm}
\usepackage{amssymb}
\usepackage{amscd}
\usepackage{graphicx}
\usepackage{epsfig}

\usepackage{latexsym, amsfonts}
\usepackage{amsthm}
\usepackage{t1enc}
\usepackage[mathscr]{eucal}
\usepackage{indentfirst}
\usepackage{graphicx, pb-diagram}
\usepackage{fancyhdr}
\usepackage{fancybox}
\usepackage{enumerate}
\usepackage[all]{xy}
\usepackage{url}
\numberwithin{equation}{section}
\newtheorem{theorem}{Theorem}[section]
\newtheorem{lemma}[theorem]{Lemma}
\newtheorem{proposition}[theorem]{Proposition}

\theoremstyle{definition}
\newtheorem{definition}[theorem]{Definition}
\newtheorem{example}[theorem] {Example}

\newtheorem{remark}[theorem]{Remark}

\numberwithin{equation}{section}

\allowdisplaybreaks
\begin{document}
\thispagestyle{empty}

\newcommand{\NN}{\ensuremath{\mathbb N}}
\newcommand{\ZZ}{\ensuremath{\mathbb Z}}
\newcommand{\CC}{\ensuremath{\mathbb C}}
\newcommand{\RR}{\ensuremath{\mathbb R}}
\newcommand{\dd}{\ensuremath{\frak{d}}}
\newcommand{\g}{\ensuremath{\frak{g}}}
\newcommand{\gl}{\ensuremath{\frak{gl}}}
\newcommand{\h}{\ensuremath{\frak{h}}}
\newcommand{\tk}{\ensuremath{\frak{k}}}
\newcommand{\p}{\ensuremath{\frak{p}}}
\newcommand{\ft}{\ensuremath{\frak{t}}}
\newcommand{\bt}{\mathbf{t}}                  
\newcommand{\bs}{\mathbf{s}}                  

\newcommand {\mcomment}[1]{{\marginpar{*}\scriptsize{\bf Marco:}\scriptsize{\ #1 \ }}}
\newcommand {\acomment}[1]{{\marginpar{*}\scriptsize{\bf Alberto:}\scriptsize{\ #1 \ }}}
\newcommand{\tJ}{\tilde{J}}
\newcommand{\Gs}{\Gamma_s}
\newcommand{\Gc}{\Gamma_c}
\newcommand{\cH}{\mathcal{H}}
\newcommand{\cL}{\mathcal{L}}
\newcommand{\cO}{\mathcal{O}}
\newcommand{\cF}{\mathcal{F}}
\newcommand{\cI}{\mathcal{I}}
\newcommand{\tP}{\tilde{P}}
\newcommand{\hP}{\hat{P}}
\newcommand{\cha}{characteristic}
\newcommand{\sh}{\sharp}
\newcommand{\un}{\underline}
\newcommand{\sctc}{\bs^{-1}(C)\cap \bt^{-1}(C)}
\newcommand{\stpttp}{\bs^{-1}(\tP)\cap \bt^{-1}(\tP)}




\markboth{A.~S.~Cattaneo and M. Zambon}{Pre-Poisson submanifolds}
$ $
\bigskip

\bigskip

\centerline{{\Large   Pre-Poisson submanifolds}}

\bigskip
\bigskip
\centerline{{\large by Alberto S. Cattaneo and Marco Zambon}}

\vspace*{.7cm}

\begin{abstract}
In this note we consider an arbitrary submanifold $C$ of a Poisson manifold $P$
and ask whether it can be embedded coisotropically in some bigger submanifold of $P$.
We define the  classes of submanifolds relevant to the question (coisotropic, Poisson-Dirac,
pre-Poisson ones), present an answer to the above question and consider the corresponding
picture at the level of Lie
groupoids, making concrete examples in which $P$ is
the dual of a Lie algebra and
   $C$ is an affine subspace. 
 
\end{abstract}

\pagestyle{myheadings}

 
\section{Introduction}\label{intro}

In this note we wish to give an analog in Poisson geometry to the following statement in
symplectic geometry. Recall that   $(P,\Omega)$ is a symplectic manifold if $\Omega$ is a closed,
non-degenerate 2 form and that a submanifold $\hat{C}$ is called coisotropic if the symplectic
orthogonal $T\hat{C}^{\Omega}$ of $T\hat{C}$ is contained in $T\hat{C}$. The statement is: if 
  $i\colon C\rightarrow P$ is any submanifold of a symplectic manifold $(P,\Omega)$,
  then there exists some symplectic submanifold $\tilde{P}$
containing $C$ as a coisotropic submanifold iff $i^*\Omega$ has constant rank.
The submanifold $\tP$ is obtained  taking any complement 
 $R\subset TP|_C$ of $TC+TC^{\Omega}$ and ``extending $C$ along $R$''. 
 Further there is a uniqueness statement ``to first order'':
 if $\tP_1$ and $\tP_2$ are as above, then there is a symplectomorphism of $P$ fixing $C$
 whose derivative at $C$ maps ${T\tP_1}|_C$ to ${T\tP_2}|_C$. This result follows using
 techniques similar to those used 
 by Marle in \cite{Ma},  and relies on a technique known as ``Moser's path method''.

The above result should not be confused with the theorem of Gotay \cite{Go}
that states the following: any presymplectic manifold (i.e. a manifold endowed with 
a constant rank closed 2-form) can be
embedded coisotropically in some symplectic manifold, which is moreover unique up to
neighborhood equivalence. The difference is that Gotay considers an abstract presymplectic 
manifold and looks for an abstract symplectic manifold in which to embed; the problem above fixes a symplectic
manifold $(P,\Omega)$ and considers only submanifolds of $P$.

In this note we ask:
\begin{itemize}
\item[1)]
 Given an arbitrary submanifold $C$ of a
Poisson manifold $(P,\Pi)$, under what conditions does there exist
some submanifold $\tilde{P}\subset P$ such that
\begin{itemize}
\item[a)] $\tilde{P}$ has a  Poisson structure induced from $\Pi$
\item[b)] $C$ is a coisotropic submanifold of $\tilde{P}$?
\end{itemize}
\item[2)]When the submanifold $\tilde{P}$ exists, is it unique up to
neighborhood equivalence (i.e. up to a Poisson diffeomorphism on
a tubular neighborhood which fixes $C$)?
\end{itemize}
 
We will see in Section \ref{embPD} that a sufficient condition is that $C$ belongs to a particular class of submanifolds
called pre-Poisson submanifolds. In that case we also have uniqueness: if
 $\tP_1$ and $\tP_2$ are as above, then there is a  Poisson diffeomorphism of 
 (a tubular neighborhood of $C$ in)
 $P$ fixing $C$
 which takes $\tP_1$  to $\tP_2$. When the Poisson structure on $P$ 
comes from a symplectic form $\Omega$, the pre-Poisson submanifolds of $P$ are 
exactly the submanifolds for which the pullback of $\Omega$ has constant rank;
hence we improve the ``uniqueness to first order'' result in the symplectic setting mentioned above
to uniqueness in a neighborhood of $C$.
 
Since the above question is essentially about when an arbitrary submanifold can be regarded as
a coisotropic one, 
 we want to motivate in Section \ref{coiso} why coisotropic submanifolds are
interesting at all. In Section \ref{pd} we will describe the submanifolds of $P$  
which inherit a Poisson structure; these are the ``candidates'' for $\tP$ as above.
Then in Section \ref{la} we will present a non-trivial example: we consider as Poisson manifold
$P$ the dual of a Lie algebra $\g$, and as submanifold $C$ either a translate
of the annihilator of a Lie subalgebra or the annihilator of some subspace of $\g$.
Finally in Section \ref{groids} we recall how to a Poisson manifold one can associate
symplectic groupoids and investigate what pre-Poisson submanifolds correspond to at the
groupoid level, discussing again the example where $P$ is the dual of a (finite
dimensional) Lie algebra. All manifolds appearing in this note are assumed to be finite dimensional.

\textbf{Acknowledgments:} 
M.Z. thanks Dirk T\"oben for discussion on symmetric pairs and related topics.
A.S.C. acknowledges partial support of SNF Grant No.~20-113439.
This work has been partially supported
by the European Union through the FP6 Marie Curie RTN ENIGMA (Contract
number MRTN-CT-2004-5652) and by the European Science Foundation
through the MISGAM program.

\section{Coisotropic submanifolds}\label{coiso}

A manifold $P$ is called \emph{Poisson manifold} if it is endowed with a bivector field
$\Pi\in\Gamma(\Lambda^2 TP)$ satisfying $[\Pi,\Pi]=0$, where $[\bullet,\bullet]$ denotes the
Schouten bracket on multivector fields. 
 Let us denote by $\sharp\colon  T^*P \rightarrow TP$ the map
given by contraction with $\Pi$. The image of $\sharp$ is a singular integrable distribution on $P$, whose
leaves are endowed with a symplectic structure that encodes the bivector field  $\Pi$. Hence 
one can think of a Poisson manifold as a manifold with a singular foliation by symplectic
leaves.

Alternatively  $P$ is a Poisson manifold if there is a Lie bracket $\{\bullet,\bullet\}$
on the space of functions satisfying the Leibniz identity\footnote{In this case one says that 
$(C^{\infty}(P),\{\bullet,\bullet\},\cdot)$ forms a
Poisson algebra.} $\{f,g\cdot h\}=\{f,g\}\cdot h+ g
\cdot \{f, h\}$.
The Poisson bracket $\{\bullet,\bullet\}$ and the bivector field $\Pi$ determine each 
other by the formula $\{f,g\}=\Pi(df,dg)$. In this note we will use
both the geometric and algebraic characterization of Poisson manifolds.

Symplectic manifolds $(P,\Omega)$ are examples of Poisson manifolds: the map $TP \rightarrow
T^*P$ given by contracting with $\Omega$ is an isomorphism, and (the negative of) its inverse
is the sharp map of the Poisson bivector field associated to $\Pi$. Connected symplectic manifolds
are exactly the Poisson manifolds whose symplectic foliation consists of just one leaf.

A second standard example, which will be used in Section \ref{la}, is the dual $\g^*$ of a   Lie algebra $\g$, as follows. A linear function $v$ on $\g^*$ can be regarded
as an element of $\g$; one defines the Poisson bracket on linear functions
as $\{v_1,v_2\}:=[v_1,v_2]$, and the bracket for arbitrary functions is determined by this in
virtue of the Leibniz rule. Duals of Lie algebras are exactly the Poisson manifolds whose
Poisson bivector field is linear. The symplectic foliation of $\g^*$ is given by the 
orbits of the coadjoint action; the origin is a symplectic leaf, and unless  $\g$ is an
abelian Lie algebra the symplectic foliation will be singular. We will discuss this example in
more detail in Section \ref{la}.\\

A submanifold $C$ of a Poisson manifold $P$ is called \emph{coisotropic} if 
$\sharp N^*C\subset TC$. Here $N^*C$ (the conormal bundle of
$C$) is defined as the annihilator of $TC$, and the singular distribution
$\sharp N^*C$ on $C$ is called the \emph{characteristic distribution}.
  Notice that if the Poisson structure
of $P$ comes from
 a symplectic form $\Omega$ then
the subbundle $\sharp N^*C$ is just the symplectic orthogonal of $TC$, so we are reduced to the
usual definition of coisotropic submanifolds in the symplectic case. If a submanifold $C$
intersects the symplectic
leaves $\cO$ of $P$ cleanly, then $C$ is coisotropic iff each intersection $C\cap \cO$ is a
coisotropic submanifold of the symplectic manifold $\cO$.
In algebraic terms we have the following characterization: a submanifold $C$ is coisotropic iff
$I_C \colon=\{f\in C^{\infty}(P)\colon f|_C=0\}$ is a Poisson subalgebra of  
$(C^{\infty}(P),\{\bullet,\bullet\},\cdot)$.
 
 In the following we want to motivate the naturality and importance of  coisotropic
 submanifolds.
 \begin{itemize}
\item Graphs of Poisson maps are coisotropic: 
\begin{proposition}[Cor. 2.2.3 of \cite{Wcoiso}]
Let $\Phi\colon(P_1,\Pi_1) \rightarrow (P_2,\Pi_2)$ be a map between Poisson manifolds. $\Phi$
is a Poisson map (i.e. $\Phi_*(\Pi_1)=\Pi_2$) iff its graph is a coisotropic submanifold of
$ (P_1\times P_2,\Pi_1-\Pi_2)$.
\end{proposition} 

\item Certain canonical quotients of coisotropic submanifolds are Poisson manifolds: 
define $F_C\colon =\{f\in  C^{\infty}(P)\colon  \{f,I_C\}\subset I_C\}$, the Poisson normalizer of $I_C$. 
It is a Poisson subalgebra of $C^{\infty}(P)$, and $I_C\subset F_C$ is a Poisson ideal.
Further notice that $F_C$ consists exactly of the functions on $P$ whose differentials annihilate the
characteristic distribution $\sharp N^*C$.
Hence we have the following statements about the quotient of $C$ by the characteristic
distribution:
\begin{proposition}
$F_C/I_C$ inherits the structure of  a Poisson algebra. Therefore
 $\un{C}\colon =C/\sharp N^*C$, if smooth, inherits the structure of a Poisson manifold so that
 $C\rightarrow \un{C}$ is a Poisson map.
\end{proposition}
Given any Poisson algebra $A$, one can ask whether it admits a deformation quantization,
i.e. if it is possible to deform the commutative multiplication ``in direction of the Poisson
bracket'' to obtain an associative product. Remarkable work of Kontsevich \cite{K}
showed that this is always possible if $A$ is the algebra of functions on a smooth 
Poisson manifold.
The Poisson algebras $F_C/I_C$ provide natural and interesting instances of Poisson algebras
which usually cannot be regarded as
algebras of functions on a smooth manifold; the problem of their deformation quantization 
has been considered in \cite{CaFeCo1,CaFeCo2}.

\item Last, a coisotropic submanifold $C$ gives rise to a Lie subalgebroid of the Lie algebroid
associated to $P$.
Recall that a \emph{Lie algebroid} is a vector bundle $E\rightarrow P$ with a Lie bracket
$[\bullet,\bullet]$ on its space of sections and a bracket preserving bundle map $\rho\colon E\rightarrow TP$
satisfying $[e_1,fe_2]=\rho(e_1)f\cdot e_2+f[e_1,e_2]$; standard examples are tangent bundles
and Lie algebras.
 Every Poisson manifold $P$ induces the
structure of a Lie algebroid on its cotangent bundle $T^*P$: the bracket is given
by $[df,dg]=d\{f,g\}$ and the bundle map $T^*P \rightarrow TP$ by $-\sharp$. We have
\begin{proposition}[Cor. 3.1.5 of \cite{Wcoiso}]
If $C\subset P$ is coisotropic then the conormal bundle $N^*C$ is a Lie subalgebroid of $T^*P$.
\end{proposition}
\end{itemize}

\section{Poisson-Dirac and cosymplectic submanifolds}\label{pd}

In virtue of the question asked in the introduction it is necessary to determine which
submanifolds $\tP$ of a Poisson manifold $(P,\Pi)$ inherit a Poisson structure. Notice that,
unlike  symplectic forms, it is usually not possible to restrict a Poisson  bivector field to a submanifold and
obtain again a bivector field. However it is possible to view a Poisson bivector field as a Dirac structure \cite{Cou}, and Dirac structures restrict to (usually not smooth) Dirac structures on
submanifolds. This point of view led to the definition below, which we phrase without reference
to Dirac structures.
 
 We first make the following remark, in which $(\cO,\Omega)$ denotes a
 symplectic leaf of $P$ and $\tP\subset P$ some submanifold: 
the linear subspace   $T_p\tP \cap T_p{\cO}$ of $(T_p{\cO},\Omega_p)$ is a symplectic subspace
iff $\sharp N_p^*\tP\cap T_p\tP=\{0\}$. In this case 
  $T\tP_p$ is endowed with a bivector field $\tilde{\Pi}_p$, obtained essentially by inverting
the non-degenerate form ${\Omega_p}|_{T_p\tP\cap T_p{\cO}}$. Now we can make sense of the following
definition (Cor.   11 of \cite{CrF}):

\begin{definition}
A submanifold $\tP$ of $P$ is called  \emph{Poisson-Dirac submanifold} if 
$\sharp N^*\tP\cap T\tP=\{0\}$ and the induced bivector field $\tilde{\Pi}$ 
on $\tP$ is a smooth.
\end{definition}
In this case the bivector field is automatically integrable (Prop. 6 of \cite{CrF}), so that $(\tP,\tilde{\Pi})$ is a Poisson
manifold. Equivalently (Def. 4 of \cite{CrF}) $\tP$ is a Poisson-Dirac submanifold if it admits
a Poisson structure for which 
the  symplectic leaves are  (connected)
intersections with the symplectic leaves $\cO$ of $P$ and so that the former are symplectic submanifolds of
the leaves $\cO$. Notice that the inclusion $\tP\rightarrow P$ is usually not a Poisson map; it is iff
$\tP$ is a   Poisson submanifold, i.e. a smooth union of symplectic leaves.

A submanifold $\tP$ satisfying $T\tP\oplus \sharp N^*\tP=TP|_{\tP}$ is called a \emph{cosymplectic
submanifold}. In this case one can show that the induced bivector field 
$\tilde{\Pi}$ on $\tP$ is automatically
smooth, hence cosymplectic submanifolds are Poisson-Dirac submanifolds.
The Poisson bracket on a cosymplectic submanifold $\tP$  is computed as follows:
$\{\tilde{f}_1,\tilde{f}_2\}_{\tP}$ is the restriction to $\tP$ of $\{f_1,f_2\}$, where the $f_i$ are extensions of $\tilde{f}_i$ to $P$
such that $df_i|_{\sharp N^*\tP}=0$.

If the Poisson structure on $P$ comes from a symplectic 2-form, then  the
Poisson-Dirac and cosymplectic submanifolds are just the symplectic submanifolds.

\section{Coisotropic embeddings in Poisson-Dirac submanifolds}\label{embPD}

Now we determine
 under what conditions on a submanifold $i\colon C\rightarrow P$  there exists
a Poisson-Dirac submanifold $\tilde{P}\subset P$  so that
  $C$ is   coisotropic in $\tilde{P}$.
We saw in the introduction that, when the Poisson structure on $P$ comes
from a  symplectic form $\Omega$, a sufficient and necessary condition 
is that  $ker(i^*\Omega)$, which  in terms of the Poisson tensor  
  is 
 $TC\cap \sharp N^*C$, has constant rank.
 In the general Poisson case however  $TC\cap \sharp N^*C$, even when it has constant rank,
 might not be a smooth distribution on $C$. In the symplectic case $ker(i^*\Omega)$
 has constant rank iff $TC+TC^{\Omega}$ has constant rank, and 
 it turns out that this is
  the right condition to 
 generalize to the Poisson case.  This motivates
\begin{definition}[Def. 2.2 of \cite{CZ}]
A submanifold $C$ of a Poisson manifold $(P,\Pi)$ is called
\emph{pre-Poisson} if the rank of $TC+\sharp N^*C$ is constant
along $C$.
\end{definition} 
Such submanifolds were first considered in \cite{CaFa1,CaFa2}. We have 
\begin{theorem}\label{emb}[Thm. 3.3 of \cite{CZ}]
Let $C$ be a pre-Poisson submanifold of a Poisson manifold
$(P,\Pi)$. Then there exists a
cosymplectic submanifold $\tP$ containing $C$ such that $C$ is
coisotropic in $\tP$.
\end{theorem}
\begin{proof}[Sketch of the proof]
Because of the rank condition on $C$ we can choose a smooth
subbundle $R$ of $TP|_C$ which is a complement  to $TC + \sharp
N^*C$.  By linear algebra, at every point $p$ of $C$,  $T_pC\oplus R_p$ is a cosymplectic subspace of $T_pP$ in
which $T_pC$ sits coisotropically. Now we ``thicken''  $C$ to a smooth
submanifold  $\tP$ of $P$ satisfying $T\tP|_C=TC\oplus R$. One can show that 
in a neighborhood of $C$
  $\tP$ is a cosymplectic submanifold, so shrinking $\tP$ if necessary we are done. \end{proof}

\begin{remark}\label{constr}
The cosymplectic submanifold $\tP$ above is constructed by 
taking any complement 
 $R\subset TP|_C$ of $TC+\sharp N^*C$ and ``extending $C$ along $R$''. 
 
There are submanifolds $C$ which are not pre-Poisson but   still admit some Poisson-Dirac
submanifold $\tP$ in which they embed coisotropically. This happens for example  
if $C$ has trivial intersection with the symplectic leaves of $P$ (and the symplectic foliation
of $P$ is not regular): in this case $\tP:=C$ is a Poisson-Dirac submanifold, the induced
 Poisson bivector field being zero.

However, if we ask that the submanifold $\tP$ be not just Poisson-Dirac but actually
cosymplectic, then $C$ is necessarily a pre-Poisson submanifold, and $\tP$ is constructed
as described above
(Lemma 4.1 of \cite{CZ}).
\end{remark}  

The following are elementary examples of pre-Poisson submanifolds and of cosymplectic
submanifolds in which they embed coisotropically. 
In section \ref{la} we will give less trivial examples; see also Section 5 of \cite{CZ}.

\begin{example}
When $C$ is a coisotropic submanifold of
$P$,
the construction of Thm. \ref{emb} delivers $\tP=P$ (or more
precisely, a tubular neighborhood of $C$ in $P$).
\end{example}
\begin{example}\label{exstrong}
When $C$ is just a point $x$ then the
construction of Thm. \ref{emb} delivers as $\tP$ any slice
through $x$ transversal to the symplectic leaf $\cO_x$.
\end{example}
\begin{example}
If $C_1\subset P_1$ and $C_2\subset P_2$ are pre-Poisson
submanifolds of Poisson manifolds, the cartesian product
$C_1\times C_2 \subset P_1\times P_2$ also is, and if the
construction of Thm. \ref{emb} gives cosymplectic submanifolds
$\tP_1\subset P_1$ and $\tP_2\subset P_2$, the same construction
applied to $C_1\times C_2$ (upon suitable choices of complementary
subbundles) delivers the cosymplectic submanifold $\tP_1\times
\tP_2$ of $P_1\times P_2$.
\end{example}
 
The following lemma will be useful in Section \ref{la}:
\begin{lemma}\label{sur}
Let $P_1,P_2$ be Poisson manifolds and 
$f\colon P_1 \rightarrow P_2$ be a submersive Poisson morphism.  
 If $C\subset P_2$ is a pre-Poisson submanifold then $f^{-1}(C)$ is a pre-Poisson submanifold of $P_1$.
Further, if $\tilde{P}_2$ is a cosymplectic submanifold containing $C$ as a coisotropic submanifold,
then $f^{-1}(\tilde{P}_2)$ is a cosymplectic submanifold containing $f^{-1}(C)$ as a coisotropic
submanifold.
\end{lemma}
\begin{proof}
Let $y\in C$ and $x\in f^{-1}(y)$. Since
$$f_*(\sharp N^*_x (f^{-1}(C)))=f_*(\sharp f^*(N^*_yC))=\sharp N^*_yC$$ it follows that
the restriction of $f_*$ to $T_x (f^{-1}(C))+\sharp N^*_x (f^{-1}(C))$ has image $T_yC+\sharp N_y^*C$,
whose rank is independent of $y \in C$ by assumption.
Since the kernel of this restriction, being $T_x(f^{-1}(y))$, also has constant rank, it follows that
$f^{-1}(C)$ is pre-Poisson. 

Further it is clear that $f_*$ maps a complement $R_x$ of $T_x (f^{-1}(C))+\sharp N^*_x (f^{-1}(C))$ in $T_x P_1$ isomorphically
onto a complement $R_y$ of $T_yC+\sharp N_y^*C$ in $T_yP_2$, so that $R_x +T_x (f^{-1}(C))$
is the pre-image of $R_y +T_y C$ under $f_*$. Using Remark \ref{constr} this proves the second assertion.
\end{proof}

The
answer to the problem of uniqueness is given by  
\begin{theorem}\label{uniPhi}[Thm. 4.4 of \cite{CZ}]
Let $C$ be a pre-Poisson submanifold $(P,\Pi)$, and $\tP_0$,
$\tP_1$ cosymplectic submanifolds that contain $C$ as a coisotropic submanifold. Then,
shrinking $\tP_0$ and $\tP_1$ to a smaller tubular
neighborhood of $C$ if necessary, there is a Poisson diffeomorphism $\Phi$ of $P$ taking
$\tP_0$ to $\tP_1$ and which is the identity on $C$.
\end{theorem}

\begin{proof}[Sketch of proof]
In a neighborhood $U$ of $\tP_0$ take a projection $\pi\colon
U\rightarrow\tP_0$.  Applying Thm. \ref{emb}
  one can construct a curve of cosymplectic submanifolds $\tP_t$ containing $C$
 which,
at points of $C$, are all transverse to  the fibers of $\pi$. Using the cosymplectic 
submanifolds $\tP_t$ one can construct a hamiltonian time-dependent vector field 
$X_{H_t}$
whose time-$t$ flow maps
$\tP_0$ to $\tP_t$. Further $X_{H_t}$ vanishes on $C$, hence its time-1 flow is the identity on $C$.
\end{proof}

\section{Duals of Lie algebras}\label{la}

In this subsection $\g$ will always denote a finite dimensional Lie algebra.
We saw in Section \ref{coiso}   that its dual $\g^*$ is a Poisson manifold, whose Poisson
bracket on linear functions (which can be identified
with elements of $\g$)
is given by $\{g_1,g_2\}:=[g_1,g_2]$.  In what follows we will need the notion of 
adjoint
action of $G$ on $\g$, which is $Ad_g v:=\frac{d}{dt}|_0 g \cdot exp(tv) \cdot g^{-1}$.
Its derivative at the identity gives the Lie algebra action of $\g$ on itself by
$ad_wv:=\frac{d}{dt}|_0 Ad_{exp(tw)}v=[w,v]$. We will also consider the (left) actions $Ad^*$ and $ad^*$ on $\g^*$
obtained by dualizing; the orbits of the coadjoint action $Ad^*$ are exactly the symplectic
leaves of the Poisson manifold $\g^*$.

It is known that if $\h$ is a Lie subalgebra of $\g$,
then its annihilator $\h^{\circ}$ is a coisotropic submanifold of $\g^*$ (also see Prop. \ref{hg} below).
We shall look at two generalizations: the first  considers affine subspaces obtained translating 
$\h^{\circ}$; the second is
obtained by weakening the condition that $\h$ be a subalgebra.\\

\begin{proposition}\label{hg}
Let $\h$ be a Lie subalgebra of $\g$ and fix $\lambda \in \g^*$. Then
the affine subspace $C:=\h^{\circ}+\lambda$ is always pre-Poisson, and it is coisotropic
iff $\lambda$ is a character of $\h$ (i.e. by definition $\lambda \in [\h,\h]^{\circ}$).
 \end{proposition}
\begin{proof}
The restriction $f\colon  \g^* \rightarrow \h^*$ is a Poisson map because $\h$ is a Lie subalgebra. Every point $\nu$ of $\h^*$ is a pre-Poisson
submanifold (see Ex. \ref{exstrong}), hence by Lemma \ref{sur} its pre-image $f^{-1}(\nu)$ (which will be a translate of 
$\h^{\circ}$) is pre-Poisson. Notice that by Lemma \ref{sur} we also know that,
for any slice $S\subset \h^{*}$ transverse to the $H$-coadjoint orbit through $\nu$, 
 $f^{-1}(S)$ is a cosymplectic submanifold containing coisotropically $f^{-1}(\nu)$.
 Further from the proof of Lemma \ref{sur} it is clear that $f^{-1}(\nu)$ is coisotropic in $\g^*$ iff $\{\nu\}$ is coisotropic in $\h^*$,
 i.e. if $\nu$ is a fixed-point of the $H$-coadjoint action or equivalently $\nu|_{[\h,\h]}=0$.
\end{proof}

\begin{example}\label{sl2}
Let $\g= \frak{sl}(2,\RR)$. In a suitable basis the Lie algebra structure is given by 
$[e_1,e_2]=-e_3,[e_2,e_3]=e_1,[e_3,e_1]=e_2$. The symplectic leaves of $\g^*$
are given essentially by the connected components of level sets of the Casimir function $\nu_1^2+
\nu_2^2-\nu_3^2$ (where $\nu_i$ is just $e_i$ viewed as a linear function on
$\g^*$), and they consist of a family of two-sheeted hyperboloids, the 
cone\footnote{The cone is the union of 3 leaves, one being
the origin.} $\nu_1^2+
\nu_2^2-\nu_3^2=0$ and a family of one-sheeted hyperboloids \cite{CW}. 
  $C:=\{(0,t,t):t\in
\RR\}\subset \g^*$ is contained in the cone and is clearly a coisotropic
submanifold; indeed it is the annihilator of the Lie subalgebra
$\h:=span\{e_1,e_2-e_3\}$ of $\g$. If we translate $C$ by an element 
in the annihilator of $[\h,\h]=\RR(e_2-e_3)$ we obtain an affine line
 contained in one of the hyperboloids, which hence is lagrangian
 there, therefore coisotropic in $\g^*$. If we translate $C$ by any other
$\lambda \in \g^*$ we obtain a line that intersects transversely the
hyperboloids, so at every point of such a line $C'$ we have $TC'+\sharp
N^*C'=T\g^*$,  showing that $C'$ is pre-Poisson.
\end{example}

Before considering the case when $\h$ is \emph{not} a subalgebra of $\g$ we need the 
\begin{lemma}\label{formula}
Let $C\subset \g^*$ be an affine subspace obtained by translating
the annihilator of a linear subspace $\h\subset \g$. Then $\sharp N^*_xC
=ad_{\h}^*(x):=\{ad_h^*(x):h\in \h\}$ for all $x\in C$.
\end{lemma}
\begin{proof}
$N^*_xC$ is given by the differentials at $x$ of the functions $h\in \h\subset
C^{\infty}(\g^*)$. Now for any $g\in \g$ we have
$$\langle \sharp d_xh, g \rangle = d_xg(\sharp d_xh)=\{h,g\}(x)=\langle [h,g], x \rangle=
\langle ad_h^*(x), g \rangle,$$ i.e. $\sharp d_xh=ad_h^*(x)$.
\end{proof}

\begin{remark}\label{alt}
An alternative proof of Prop. \ref{hg} can be given using Lemma \ref{formula}. Indeed
any $x\in C$ can be written uniquely as $y+\lambda$ where $y\in \h^{\circ}$.
Notice that $ad_h^*(y)\in  \h^{\circ}$ for all $h\in \h$, because $\langle ad_h^*(y), \h
\rangle = \langle y, [h,\h] \rangle$ vanishes since $\h$ is a subalgebra. Hence
$$T_xC+\sharp N^*_x C= \h^{\circ}+ \{ad_h^*(y)+ad_h^*(\lambda):h\in \h\}=
\h^{\circ}+  ad_{\h}^*(\lambda),$$ which is independent on the point $x$. 
From the first computation above (applied to $\lambda$ instead of  $y$)
 it is clear that $ad_{\h}^*(\lambda)\in \h^{\circ}$
iff $\lambda \in [\h,\h]^{\circ}$.
\end{remark}

Now we consider the case when $\h$ is just a  linear subspace of $\g$
and $\h^{\circ}\subset \g^*$ its
dual. Since the Poisson tensor of $\g^*$ vanishes at the origin we have
$T(\h^{\circ})+\sharp N^*(\h^{\circ})=T(\h^{\circ})$ at the origin, so $\h^{\circ}$ is pre-Poisson iff it is coisotropic
(i.e. if $\h$ is a Lie subalgebra). 
The  open subset $C$ of $\h$ on which
$T(\h^{\circ})+\sharp N^*(\h^{\circ})$ has maximal rank will be pre-Poisson.
Then, shrinking $C$ if necessary,  
we can find a subspace $R\subset \g^*$ (independent of $x\in C$) with $R\oplus(T_xC+\sharp
N_x^*C)=\g^*$ for all $x\in C$. For example 
we can construct such an $R$ at one point $\bar{x}$ of $C$, and  since transversality is an open condition,
$R$ will be transverse to $TC+\sharp
N^*C$ in a neighborhood of $\bar{x}$ in $C$.
By Thm. \ref{emb} an open
 subset $\tilde{P}$  of the subspace $\p^{\circ}:=R\oplus C$ (containing $C$) is cosymplectic. 
 If we
\emph{assume} that    
$\sharp N_y^*\tilde{P}$ is independent of the point $y\in \tilde{P}$  we are in the situation of
the following proposition. 

\begin{proposition}\label{dual}
Let $\p$ be a linear subspace of $\g$ such that  an open subset $\tilde{P}\subset\p^{\circ}$
is cosymplectic and $\tk^{\circ}:=\sharp N_y^*\tilde{P}$ is independent of   $y\in \tilde{P}$.
 Then $\tk\oplus \p=\g$, $\tk$ is a Lie subalgebra of $\g$
and $[\tk,\p]\subset \p$. Hence, whenever $[\p,\p]\subset \tk$, 
$(\tk,\p)$ forms a symmetric pair \cite{Kn}. 
\end{proposition}
\begin{proof}
The fact that $\tk\oplus \p=\g$ follows from  $\tk^{\circ}\oplus \p^{\circ}=\g^*$,
which holds because $\tilde{P}$ is cosymplectic.
Recall that given    functions $f_1,f_2$ on $\tilde{P}$, the bracket
$\{f_1,f_2\}_{\tilde{P}}$ is obtained by extending the functions in a constant 
way along
$\tk^{\circ}$ to obtain functions $\hat{f}_1, \hat{f}_2$ 
on $\g^*$, taking their Poisson bracket
and restricting to $\tilde{P}$. 
 Further (see Cor. 2.11 of \cite{Xu}) the differential of 
$\{\hat{f}_1, \hat{f}_2\}$ at any point of $\tilde{P}$ 
  kills $\tk^{\circ}$. 
So if the $f_i$ are restrictions of linear functions on
$\p^{\circ}$ then  $\hat{f}_i$ will 
be linear functions on $\g^*$ 
corresponding to elements of $\tk$, 
and
$\{\hat{f}_1, \hat{f}_2\}$, which is a linear function on $\g^*$, will also
correspond to an element of $\tk$.
 We deduce that $\tk$ is a Lie subalgebra of $\g$
(and that the    Poisson structure on $\tilde{P}$ induced from $\g^*$ 
is the restriction of a linear Poisson
structure on  $\p^{\circ}$). 

To show $[\tk,\p]\subset \p$ pick any $k\in \tk ,p\in \p$ and $y\in
\tilde{P}$. Then $\langle [k,p], y \rangle=
-\langle k, ad^*_p(y) \rangle=
  \langle k,\sharp d_yp  \rangle=0,$
  using Lemma \ref{formula} in the second equality,
because $\sharp d_yp \subset \sharp N_y^*\tilde{P}=\tk^{\circ}$.
This shows that $[k,p]$ annihilates $\tilde{P}$, hence it must annihilate
its span $\p^{\circ}$.
\end{proof}
\begin{remark}
The text preceding Prop. \ref{dual} and the proposition itself 
give a way to start with a simple piece of data (a subspace of $\g$) and,
in favorable cases, obtain a decomposition 
$\tk\oplus \p=\g$ where $\tk$ is a Lie subalgebra  
and $[\tk,\p]\subset \p$. 
If $\g$ admits a non-degenerate $Ad$-invariant bilinear form $B$,
 then the $B$-orthogonal $\p$ of
any subalgebra $\tk$ satisfies $[\tk,\p]\subset \p$, because for any
$k,k'\in \tk$ and $p\in \p$ we have
$B([k,p],k')=-B(p,[k,k'])=0$. 
If $B$ is positive-definite (such a $B$ exists for example if  
the simply connected Lie group integrating $\g$ is compact), then we clearly also have 
$\tk\oplus \p=\g$. Hence for such Lie algebras one obtains the kind of decomposition of Prop. \ref{dual} in a much easier way.
\end{remark}

A converse statement to Prop. \ref{dual} is given by 
 \begin{proposition}
Assume that $\tk\oplus \p=\g$,
$[\tk,\p]\subset \p$ and  there exists a point $y\in  \p^{\circ}$ at which none of
the fundamental vector fields $\frac{d}{dt}|_0 Ad^*_{exp(tp)}(y)$  
 vanish, where $p$ ranges over $\p \setminus \{0\}$. Then 
there is  an open subset $\tilde{P}\subset\p^{\circ}$
which is cosymplectic and $\tk^{\circ}:=\sharp N_x^*\tilde{P}$ is independent of   $x\in \tilde{P}$.
(Hence applying Prop. \ref{dual} it follows that $\tk$ is a Lie subalgebra of $\g$).
 \end{proposition}
 \begin{proof} 
 For all $x\in \p^{\circ}$ we have $\sharp N^*_x (\p^{\circ})=ad_{\p}^*(x) \subset \tk^{\circ}$, as can be seen using
    $\langle ad_p^*(x), \tk \rangle =
 \langle  x, [p,\tk] \rangle =0$ for all $p\in \p$ (which holds because of $[\tk,\p]\subset \p$).
 The assumption on the coadjoint action at $y$ means that the map
$\p \rightarrow \g^*, p \mapsto ad^*_p(y)$ is injective; by continuity it is injective also on an 
open subset $\tilde{P}\subset\p^{\circ}$, 
 and by dimension counting 
we get $\sharp N^*_{x} (\p^{\circ})=\tk^{\circ}$ on $\tilde{P}$.
\end{proof}

Now we display an example for Prop. \ref{dual}
\begin{example}
Let $\g=\gl(2,\RR)$. We identify $\g$ with $\g^*$ via the non-degenerate (indefinite) inner
product $(A,B)=Tr(A\cdot B)$. Since it is   $Ad$-invariant, the action of $ad_X$ and $ad^*_X$ on 
$\g$ and $\g^*$ are intertwined (up to sign).

Now take $\h=\left \{ \left( \begin{smallmatrix}  0 & b \\ c & d
\end{smallmatrix} \right):b,c,d \in \RR \right \}$, which is not a subalgebra. Its annihilator   is identified with the line $C$ spanned by 
 $\left(\begin{smallmatrix}  1 & 0 \\ 0 & 0
\end{smallmatrix}\right)$. Since $C$ is one-dimensional and the Poisson structure on $\g^*$ linear it is 
clear that  $\sharp N^*_x C$ is independent of $x\in C\setminus\{0\}$ and
$C\setminus\{0\}$ is pre-Poisson.
Using Lemma \ref{formula} we compute $\sharp N^*_x C=
\left \{\left(\begin{smallmatrix}  0 & b \\ c & 0
\end{smallmatrix}\right):b,c \in \RR \right \}$, so as complement $R$ to $T_xC+\sharp N^*_x C$ we can 
take the line spanned
by $\left(\begin{smallmatrix}  0 & 0 \\ 0 & 1
\end{smallmatrix}\right)$. Then $\p^{\circ}:=R\oplus C$ is given by the diagonal matrices, and $\p\subset \g$
is given by matrices with only zeros on the diagonal.
For any $\left(\begin{smallmatrix}  a & 0 \\ 0 & d
\end{smallmatrix}\right)\in \p^{\circ}$ we compute $\sharp N^*_{\left(\begin{smallmatrix}  a & 0 \\ 0 & d
\end{smallmatrix}\right)} \p^{\circ}$ using Lemma \ref{formula} and obtain the set of matrices with only 
zeros on the
diagonal if $a\neq d$ and $\{0\}$ otherwise.
So the open set $\tilde{P}$ on which $\p^{\circ}$ is cosymplectic is  
a plane with a  line removed, and
$\tk^{\circ}:=\sharp N_{\left(\begin{smallmatrix}  a & 0 \\ 0 & d
\end{smallmatrix}\right)}^*\tilde{P}$ is independent of the footpoint $\left(\begin{smallmatrix}  a & 0 \\
0 & d
\end{smallmatrix}\right)\in \tilde{P}$. 
$\tk\subset \g$ coincides hence with the set of diagonal matrices. 
As predicted by Lemma \ref{dual}   $\tk$ is a Lie subalgebra  
and $[\tk,\p]\subset \p$; one can check easily that $[\p,\p]\subset \tk$ too.

Since $\tk$ is abelian, the linear Poisson structure induced on $\tilde{P}$   is the zero Poisson structure. This can be seen also looking at the explicit Poisson
structure on $\g^*$, which with respect to the coordinates  given by the basis $a=
\left(\begin{smallmatrix}  1 & 0 \\ 0 & 0
\end{smallmatrix}\right)$, $b=
\left(\begin{smallmatrix}  0 & 1 \\ 0 & 0
\end{smallmatrix}\right)$, $c=
\left(\begin{smallmatrix}  0 & 0 \\ 1 & 0
\end{smallmatrix}\right)$ and 
$d=
\left(\begin{smallmatrix}  0 & 0 \\ 0 & 1
\end{smallmatrix}\right)$ of $\g^*$ is   
$$ 
-b\partial_a\wedge \partial_b+
c\partial_a\wedge \partial_c +
(d-a)\partial_b\wedge \partial_c
-b\partial_b\wedge \partial_d
+c\partial_c\wedge \partial_d.$$
Indeed at a point $\left(\begin{smallmatrix}  a & 0 \\ 0 & d
\end{smallmatrix}\right)$ of  $\p^{\circ}$ the bivector field reduces to $(d-a)\partial_b\wedge \partial_c$. 
Finally remark that if we had chosen $R$ to be spanned by $\left(\begin{smallmatrix}  0 & 0 \\ 1 & 1
\end{smallmatrix}\right)$ instead we would have obtained as
$\sharp N^*_{\left(\begin{smallmatrix}  a & b \\ 0 & b
\end{smallmatrix}\right)}\p^{\circ}$ the span of $\left(\begin{smallmatrix}  -b & b \\ a-b & b
\end{smallmatrix}\right)$ and $\left(\begin{smallmatrix}  0 & b-a \\ 0 & 0
\end{smallmatrix}\right)$, which obviously is not constant on any open subset of $\p^{\circ}$
 \end{example}

\section{Subgroupoids associated to pre-Poisson submanifolds}\label{groids}

In Section \ref{coiso} we defined Lie algebroids and recalled that for every Poisson
manifold $P$ there is an associated Lie algebroid, namely the cotangent bundle $T^*P$.

In analogy to the fact that finite dimensional Lie algebras integrate to  Lie groups (uniquely if
required to be simply connected), Lie algebroids - when integrable - 
integrate to objects called \emph{Lie groupoids}. Recall  that a Lie groupoid over $P$
 is given by a manifold $\Gamma$ endowed with surjective submersions $\bs$,$\bt$ (called source and
 target) to the base manifold $P$, a smooth associative multiplication 
 defined on elements $g,h\in \Gamma$ satisfying $\bs(g)=\bt(h)$, an embedding
 of $P$ into $\Gamma$ as the spaces of ``identities'' and a smooth inversion map $\Gamma
 \rightarrow
 \Gamma$; see for example   \cite{MW} for the precise definition.
 The total space of the Lie algebroid associated to the Lie groupoid $\Gamma$ is
$ker(\bt_*|_{P})\subset T\Gamma|_P$, with a bracket on sections defined using left invariant
vector fields on $\Gamma$ and $\bs_*|_{P}$ as anchor. A Lie algebroid $A$ is said to be integrable if there exists a Lie
groupoid whose associated Lie algebroid is isomorphic to $A$; in this case there is a unique (up to isomorphism)
 integrating Lie
groupoid with simply connected source fibers.
 
The cotangent bundle $T^*P$ of a Poisson manifold $P$ carries more data then just a Lie
algebroid structure; when it is integrable, the corresponding Lie groupoid $\Gamma$
is actually a \emph{symplectic groupoid} \cite{MX}, i.e. \cite{MW} there is a symplectic 
form $\Omega$ on $\Gamma$ such that the graph of the multiplication is a lagrangian submanifold
of $(\Gamma \times\Gamma \times\Gamma , \Omega\times\Omega\times(-\Omega))$. $\Omega$
is uniquely determined (up to symplectic groupoid automorphism)
by the requirement that $\bt\colon  \Gamma\rightarrow P$ be a Poisson map; a canonical Lie algebroid isomorphism between $T^*P$ and $ker(\bt_*|_{P})$ is then given by mapping
$du$ (for $u$  a function on $P$) to 
the hamiltonian vector field $-X_{\bs^*u}$.
 For example,
if $P$ carries the zero Poisson structure, then the symplectic groupoid is $T^*P$
with   the canonical symplectic structure and fiberwise addition as multiplication.
We will describe in Example \ref{isot} below the symplectic groupoid of the dual of a Lie
algebra.\\

 In this Section we want to investigate how a pre-Poisson submanifold 
  $C$ of  a Poisson manifold
$(P,\Pi)$ gives rise to subgroupoids of the source simply connected symplectic groupoid $\Gamma$ 
(assuming that $T^*P$ is an integrable Lie algebroid).
By Prop. 3.6 of \cite{CZ} $N^*C\cap
\sharp^{-1}TC$ is a Lie subalgebroid of $T^*P$. When $\sharp N^*C$
has constant rank there is another Lie subalgebroid
associated to $C$, namely $\sharp^{-1} TC=(\sharp N^*C)^{\circ}$.   We
want to describe the  subgroupoids\footnote{Here, for any Lie
subalgebroid $A$ of $T^*P$ integrating to a source simply connected Lie groupoid
$H$, we take ``subgroupoid'' to mean the (usually just immersed)
image of the (usually not injective) morphism $H\rightarrow
\Gamma$ induced from the inclusion $A\rightarrow T^*P$.}
 of
 $\Gamma$ integrating $N^*C\cap
\sharp^{-1}TC$ and $\sharp^{-1} TC$. 
 
\begin{proposition}\label{groid}[Prop. 7.2 of \cite{CZ}]
Let $C$ be a pre-Poisson submanifold of $(P,\Pi)$. Then the
subgroupoid of $(\Gamma,\Omega)$ integrating $N^*C\cap \sharp^{-1}TC$ is an
isotropic subgroupoid.
\end{proposition}

We exemplify Prop. \ref{groid} 
 considering the dual of a Lie algebra $\g$ as a Poisson manifold, as in Section \ref{la}. 
The symplectic groupoid of $\g^*$  (see Ex. 3.1 of \cite{MW})
is $T^*G$ with its canonical symplectic form, where $G$ is the simply connected Lie group
integrating $\g$. To describe the groupoid structure we identify $T^*G$ with $\g^*\times G$
by (the cotangent lift of) right translation. Then the target map $\g^*\times G \rightarrow
\g^*$ is $\bt(\xi,g)=\xi$ and the source map is $\bs(\xi,g)=Ad^*_{g^{-1}} \xi$,
and the multiplication is $(\xi,g_1)\cdot (Ad^*_{g^{-1}} \xi,g_2)=(\xi, g_1g_2)$.

\begin{example}\label{isot}
Let $\h$ be a Lie subalgebra of $\g$ and $\lambda \in \g^*$.  By Prop. \ref{hg} we know that
$C:=\h^{\circ}+\lambda$ is a pre-Poisson submanifold of $\g^*$. We claim here that the
subgroupoid of $\g^*\times G$  integrating the Lie subalgebroid $N^*C\cap \sharp^{-1}TC$
is $C\times D$, where  the subgroup $D\subset G$ is the connected component of the identity of 
$\{g\in H: (Ad^*_g \lambda)|_{\h}=\lambda|_{\h}\}$. By Prop. \ref{groid} we know that it is 
an isotropic subgroupoid.

To prove our claim, we first make the Lie subalgebroid more explicit: for all $x\in C$ 
using Remark \ref{alt} 
we have
 $$N_x^*C\cap \sharp^{-1}T_xC=
 ( \h^{\circ}+  ad_{\h}^*(\lambda))^{\circ}=\h \cap \{v\in \g: (ad^*_v \lambda)|_{\h}=0\} =:\dd,$$
 so that the Lie subalgebroid $N^*C\cap \sharp^{-1}TC\subset T^*\g^*=\g^*\times \g$
 is just the product $C\times \dd$.
 The canonical Lie algebroid isomorphism  $T^*P\cong ker(\bt_*|_{P}), 
 du\mapsto -X_{\bs^*u}$ is just the identity on $\g^* \times \g$, as can be checked
 using the explicit formula for the symplectic form on the  groupoid $\g^*\times G$
 given in Ex. 3.1 of \cite{MW}. Now
  notice that the Lie subalgebra $\dd$ integrates to
 the connected subgroup $D$ defined above. Using the definition of $D$ one
  checks that $\bt$ and $\bs$ map  $C\times D$ into $C$, and the fact that $D$ is a subgroup
  allows us to check
 that  $C\times D$ is actually 
  a Lie subgroupoid of $\g^*\times G$, proving our claim. 
\end{example}

Now we consider $\sharp^{-1}TC$ and assume that it has constant rank,
 or equivalently that
 the
characteristic distribution $TC\cap \sh N^*C$ have constant
rank\footnote{Indeed more generally we have the following for any
submanifold $C$ of $P$: if any two of $\sharp^{-1}TC$, $\sharp
N^*C+TC$ or $TC\cap \sh N^*C $ have constant rank, then the
remaining one also has constant rank. This follows trivially from
$rk (\sharp N^*C+TC)= rk (\sh N^*C)+ \dim C - rk (TC\cap \sh N^*C)
$.}.
Then $\sharp^{-1}TC$ is a Lie subalgebroid of $T^*P$, and 
quoting part of Prop. 7.2 of \cite{CZ}:
\begin{proposition}\label{pre}
The subgroupoid of $\Gamma$ integrating  $\sharp^{-1}TC$
is
 $\sctc$, and it is  a  
presymplectic  submanifold  of $(\Gamma,\Omega)$.
\end{proposition}

 \begin{remark}
 In this case the foliation integrating the characteristic distribution of $\sctc$ (i.e. the kernel of the pullback of $\Omega$) is given by the orbits of the action by right and left multiplication of
 the source-connected isotropic subgroupoid integrating $N^*C\cap \sharp^{-1}TC$.
 \end{remark}

 \begin{example}
Let $C$ be a submanifold of $\g^*$ such that $T_xC \cap T_x \cO=\{0\}$ 
at every point $x$ where $C$ intersects a 
coadjoint orbit $\cO$.  Then $C$ is pre-Poisson iff $\sharp^{-1}TC$
has constant rank, which in this case just means that the coadjoint orbits that
$C$ intersects all have the same dimension. By the above proposition the source connected
subgroupoid of $\g^*\times G$ integrating $\sharp^{-1} TC$ is 
$\{(x,g): x\in C, Ad^*_{g^{-1}}(x)=x\}$,   a bundle of groups integrating a bundle of isotropy Lie algebras of
the coadjoint action. We also have the following alternative description for this bundle of Lie algebras, which sometimes is more convenient for computations:
 $\sharp^{-1}T_xC=(\sharp N_x^*C)^{\circ}$
  can be described as $N_x^*\cO$, for $\cO$ the coadjoint
orbit through $x$.
  
If $\h$ is a Lie subalgebra of $\g$ and $\lambda \in \g^*$, we know that $C:=\h^{\circ}+\lambda$ is a pre-Poisson submanifold of $\g^*$, but
generally $\sharp^{-1}TC $ does not have constant rank. A case where it has a constant rank is the following.
As in Example
\ref{sl2} consider
 $\g= \frak{sl}(2,\RR)$ and  the pre-Poisson submanifold $C:=\{(0,t,t+1):t\in
\RR\}$. As remarked there $C$ intersects transversely the symplectic leaves of $\g^*$,
which are the  level sets of the 
 Casimir function $\nu_1^2+
\nu_2^2-\nu_3^2$. At $x=(0,t,t+1)$ we have $N^*_x\cO=\RR( t d \nu_2-(t+1) d \nu_3)$,
which in terms of the basis $e_1=\frac{1}{2}\left( \begin{smallmatrix}  0 & 1 \\ 1 & 0
\end{smallmatrix} \right)$, $e_2 =\frac{1}{2}\left( \begin{smallmatrix}  1 & 0 \\ 0 & -1
\end{smallmatrix} \right)$ and $e_3=\frac{1}{2}\left( \begin{smallmatrix}  0 & 1 \\ -1 & 0
\end{smallmatrix} \right)$ of $\frak{sl}(2,\RR)$ used in Example \ref{sl2} is
 $\RR\left( \begin{smallmatrix}  t & -(t+1) \\ t+1 & -t \end{smallmatrix} \right)$.
  As seen above, integrating these Lie algebras to subgroups of $G$ (the 
 simply connected Lie group integrating $\frak{sl}(2,\RR)$) 
  we obtain the  source connected
subgroupoid of $\g^*\times G$ integrating $\sharp^{-1}TC$.
\end{example}

\bibliographystyle{habbrv}
\bibliography{bibjapan}

\begin{thebibliography}{10}

\bibitem{CaFa1}
I.~Calvo and F.~Falceto.
\newblock Poisson reduction and branes in {P}oisson-sigma models.
\newblock {\em Lett. Math. Phys.}, 70(3):231--247, 2004.

\bibitem{CaFa2}
I.~Calvo and F.~Falceto.
\newblock {Star products and branes in Poisson-Sigma models}.
\newblock {\em Commun. Math. Phys.}, 268(3):607--620, 2006.

\bibitem{CW}
A.~Cannas~da Silva and A.~Weinstein.
\newblock {\em Geometric models for noncommutative algebras}, volume~10 of {\em
  Berkeley Mathematics Lecture Notes}.
\newblock American Mathematical Society, Providence, RI, 1999.

\bibitem{CaFeCo1}
A.~S. Cattaneo and G.~Felder.
\newblock Coisotropic submanifolds in {P}oisson geometry and branes in the
  {P}oisson sigma model.
\newblock {\em Lett. Math. Phys.}, 69:157--175, 2004.

\bibitem{CaFeCo2}
A.~S. Cattaneo and G.~Felder.
\newblock {Relative formality theorem and quantisation of coisotropic
  submanifolds}.
\newblock {\em Adv. in Math.}, 208:521--548, 2007.

\bibitem{CZ}
A.~S. Cattaneo and M.~Zambon.
\newblock {Coisotropic embeddings in Poisson manifolds}, to appear in
  \emph{Trans. Amer. Math. Soc.}

\bibitem{Cou}
T.~J. Courant.
\newblock Dirac manifolds.
\newblock {\em Trans. Amer. Math. Soc.}, 319(2):631--661, 1990.

\bibitem{CrF}
M.~Crainic and R.~L. Fernandes.
\newblock Integrability of {P}oisson brackets.
\newblock {\em J. Differential Geom.}, 66(1):71--137, 2004.

\bibitem{Go}
M.~J. Gotay.
\newblock On coisotropic imbeddings of presymplectic manifolds.
\newblock {\em Proc. Amer. Math. Soc.}, 84(1):111--114, 1982.

\bibitem{Kn}
A.~W. Knapp.
\newblock {\em Lie groups beyond an introduction}, volume 140 of {\em Progress
  in Mathematics}.
\newblock Birkh\"auser Boston Inc., Boston, MA, second edition, 2002.

\bibitem{K}
M.~Kontsevich.
\newblock Deformation quantization of {P}oisson manifolds.
\newblock {\em Lett. Math. Phys.}, 66(3):157--216, 2003.

\bibitem{MX}
K.~C.~H. Mackenzie and P.~Xu.
\newblock Integration of {L}ie bialgebroids.
\newblock {\em Topology}, 39(3):445--467, 2000.

\bibitem{Ma}
C.-M. Marle.
\newblock Sous-vari\'et\'es de rang constant d'une vari\'et\'e symplectique.
\newblock In {\em Third Schnepfenried geometry conference, Vol. 1
  (Schnepfenried, 1982)}, volume 107 of {\em Ast\'erisque}, pages 69--86. Soc.
  Math. France, Paris, 1983.

\bibitem{MW}
K.~Mikami and A.~Weinstein.
\newblock Moments and reduction for symplectic groupoids.
\newblock {\em Publ. Res. Inst. Math. Sci.}, 24(1):121--140, 1988.

\bibitem{Wcoiso}
A.~Weinstein.
\newblock Coisotropic calculus and {P}oisson groupoids.
\newblock {\em J. Math. Soc. Japan}, 40(4):705--727, 1988.

\bibitem{Xu}
P.~Xu.
\newblock Dirac submanifolds and {P}oisson involutions.
\newblock {\em Ann. Sci. \'Ecole Norm. Sup. (4)}, 36(3):403--430, 2003.

\end{thebibliography}



\noindent 

A.S. Cattaneo and M. Zambon\\
Institut f\"ur Mathematik, Universit\"at Z\"urich-Irchel, Winterthurerstr. 190, 
CH-8057 Z\"urich, Switzerland\\
\texttt{alberto.cattaneo@math.unizh.ch, marco.zambon@math.unizh.ch.}
        

\label{lastpage}
\end{document}